\documentclass[11pt]{article}
\usepackage{amssymb,amsmath,amsfonts,dsfont,mathrsfs,euscript,eufrak,colonequals,indentfirst}
\usepackage{textcomp}
\usepackage{latexsym}

\usepackage{authblk}
\newtheorem{Theorem}{Theorem}
\newtheorem{Lemma}{Lemma}
\newtheorem{Proposition}{Proposition}

\newcommand{\N}{\mathbb N}
\newcommand{\R}{\mathbb R}

\newcommand{\Expec}{\mathrm{\mathbf{E}}}

\newcommand{\defeq}{\colonequals}

\newcommand{\dd}{\mathrm{d}}
\newcommand{\id}{\mathds{1}}

\newcommand{\CovFunc}{\mathcal{K}}

\newcommand{\Probab}{\mathrm{\mathbf{P}}}
\newcommand{\Var}{\mathrm{\mathbf{Var}}}

\newcommand{\e}{\varepsilon}

\begin{document}
\title{Asymptotic analysis in multivariate worst case approximation with Gaussian kernels}
\author{A. A. Khartov$^{1,2,3}$, I. A. Limar$^{2,3}$}

\footnotetext[1]{Laboratory for Approximation Problems of Probability, Smolensk State University, 4 Przhevalsky st., 214000 Smolensk, Russia. }
\footnotetext[2]{Saint-Petersburg National Research University of Information Technologies, Mechanics and Optics (ITMO University), 49 Kronverksky Pr., 197101 Saint-Petersburg, Russia.}
\footnotetext[3]{Email addresses: \texttt{alexeykhartov@gmail.com} (A. A. Khartov), \texttt{ivan.limar95@gmail.com} (I. A. Limar)}

\maketitle


\begin{abstract}
	
	We consider a problem of approximation of $d$-variate functions defined on $\R^d$ which belong to the Hilbert space with tensor product-type reproducing Gaussian kernel with constant shape parameter. Within worst case setting, we investigate the growth of the information complexity as $d\to\infty$. The asymptotics are obtained for the case of fixed error threshold and for the case when it goes to zero as $d\to\infty$.  
\end{abstract}


\section{Introduction and problem setting}
We consider multivariate approximation problems
in the worst case setting for functions of $d$ variables from the Hilbert space with the squared-exponential
reproducing kernel as  $d \to \infty$, when the error threshold $\e\in(0,1)$ is   fixed or $\e=\e_d\to 0$. 

For every $d\in\N$ (set of positive integers) we consider the Hilbert space $H_{d,\gamma}$  with the following squared-exponential reproducing kernel
\begin{eqnarray}
    \label{kernel-eqnarray}
    \CovFunc_d(t, s) = \prod_{j = 1}^d \exp \left\{-\gamma^2(t_j - s_j)^2\right\},
\end{eqnarray}
where $t=(t_1,\ldots, t_d)$ and $s=(s_1, \ldots, s_d)$  are from $\R^d$ ($\R$ is the set of real numbers), $\gamma > 0$ is a shape parameter. The space $H_{d,\gamma}$ is well studied and it is widely used in numerical computations,  statistical leaning and
engineering, see \cite{ForrSobKean}, \cite{HastieTibshFried}, \cite{RasWill}, \cite{SchSmola}, \cite{Wendland}, \cite{Zhu}.

We consider the multivariate approximation problem $\mathrm{APP}_d: H_{d,\gamma} \to L_{2,d}$ with
$\mathrm{APP}_d f = f$ for all $f \in H_{d,\gamma}$, where $L_{2,d}$ is the space of  functions that have the finite norm
\begin{eqnarray*}
    \|f\|_{L_{2,d}}\defeq \Biggl(\int_{\R^d} f^2(x) \prod\limits_{j=1}^{d} \dfrac{e^{-x_j^2}}{\sqrt{\pi}}\dd x\Biggr)^{1/2},
\end{eqnarray*}
which the space is equipped with. We set $x=(x_1,\ldots, x_d)$ in the integral.
We approximate $\mathrm{APP}_d f$ by algorithms that use finitely many, say $n\in\N$, values of linear functionals. It is well known, see \cite{NovWoz1}, that we can restrict ourselves to $n$-rank linear algorithms from the following class 
\begin{eqnarray*}
    \mathcal{A}_{d, n} = \left\{\sum_{k = 1}^n a_k\, \ell_k f: a_k \in L_{2,d}, \ell_k \in H_{d,\gamma}^*\right\}.
\end{eqnarray*}
The worst case error of the algorithm $A_{d, n} \in \mathcal{A}_{d, n}$ is defined as
\begin{eqnarray*}
    e(A_{d, n}) \defeq\sup\limits_{\Vert f \Vert_{H_{d,\gamma}} \leqslant 1} \Vert \mathrm{APP}_d f - A_{d, n} f \Vert_{L_{2,d}},
\end{eqnarray*}
where $\|\,\cdot\,\|_{H_{d,\gamma}}$ is the norm of $H_{d,\gamma}$.
The $n^{th}$ minimal worst case error has the following form
\begin{eqnarray*}
    e(d, n) = \inf_{A_{d, n} \in \mathcal{A}_{d, n}}e(A_{d, n}).
\end{eqnarray*}
We also deal with the initial error, i.e we consider the value
\begin{eqnarray*}
    e(d, 0) = \sup\limits_{\Vert f \Vert_{H_{d,\gamma}} \leqslant 1}\Vert f \Vert_{L_{2,d}} = \Vert \mathrm{APP}_d\Vert,
\end{eqnarray*}
which can be thought as the error of the identical zero algorithm. We define the \textit{worst case information complexity} $n(d, \e)$ for normalized criterion (\textit{information complexity} for short) as
\begin{eqnarray}
    \label{approx-compl-def-eqnarray}
    n(d, \e) = \min\bigl\{n \in \N: e(d, n) \leqslant \e e(d, 0)\bigr\},\quad \e \in (0, 1),\quad d\in\N.
\end{eqnarray}
The multivariate worst case setting and its modifications are comprehensively described in \cite{NovWoz1, NovWoz2, NovWoz3}.

The quantity $n(d, \e)$ admits a representation via eigenvalues of the integral operator $I_d: H_{d,\gamma} \to H_{d,\gamma}$ that acts as follows
\begin{eqnarray*}
    (I_d f)(t) =    \int_{\R^d} f(s) \prod_{j = 1}^d \exp \left\{-\gamma^2(t_j - s_j)^2\right\} \prod\limits_{j=1}^{d} \dfrac{e^{-s_j^2}}{\sqrt{\pi}}\,\dd s, \quad f \in H_{d,\gamma},
\end{eqnarray*}
where $t=(t_1,\ldots, t_d)$ and $s=(s_1, \ldots, s_d)$  are from $\R^d$.
Let   $(\lambda_{d, m})_{m \in \N}$ be the  sequence of  eigenvalues of $I_d$ that is assumed to be ranked in the non-increasing order. Let $(\psi_{d, m})_{m \in \N}$ be the corresponding sequence of
orthonormal eigenfunctions of $I_d$, i.e. $(I_d \psi_{d, m})(t)=\lambda_{d, m}\psi_{d, m}(t)$, $m\in\N$. It is well known  (see \cite{NovWoz1}) that the following $n$-rank algorithm $S_{d, n}$ 
\begin{eqnarray*}
    S_{d, n}f = \sum_{j = 1}^n (f, \psi_{d, m})_{H_{d,\gamma}} \psi_{d, m},\quad f \in H_{d,\gamma},
\end{eqnarray*}
has the minimal worst case error $e(d, n)=\sqrt{\lambda_{d, n + 1}}$. The initial error is also known:  $e(d, 0) = \sqrt{\lambda_{d, 1}}$. Thus, representation \eqref{approx-compl-def-eqnarray} is written as the follows
\begin{eqnarray}
    \label{approx-compl-repr1-eqnarray}
    n(d, \e) = \min\bigl\{n \in \N: \lambda_{d, n + 1} \leqslant \e^2 \lambda_{d, 1}\bigr\},\quad \e \in (0, 1), \quad d \in \N.
\end{eqnarray}

It is a nice fact that $(\lambda_{d, m})_{m\in\N}$ and $(\psi_{d, m})_{m \in \N}$ are known for $I_d$ with kernel \eqref{kernel-eqnarray}. For $d=1$ we introduce the following convenient notation: $\lambda_k\defeq \lambda_{1,k}$ and $\psi_k\defeq \psi_{1,k}$, $k\in\N$. So for this case we have (see  \cite{NovWoz3} p. 17 and \cite{RasWill} p. 97):
\begin{eqnarray}\label{onedim-eigenvalue-eqnarray}
\lambda_k = (1 - \omega)\,\omega^{k - 1}, \quad    k\in\N,\quad \text{with}\quad	\omega \defeq \dfrac{2\gamma^2}{1+2\gamma^2+\sqrt{1+4\gamma^2}},
\end{eqnarray}  
and
\begin{eqnarray*}
	\psi_{k}(t)\defeq \sqrt{\tfrac{(1+4\gamma^2)^{1/4}}{2^{k-1} (k-1)!}}\cdot \exp\Bigl\{-\tfrac{2\gamma^2t^2}{1+\sqrt{1+4\gamma^2}}\, \Bigr\}\cdot H_{k-1}\bigl((1+4\gamma^2)^{1/4}t\bigr),\quad t\in \R,\quad k\in\N,
\end{eqnarray*}
where $H_{k-1}$ is the standard Chebyshev--Hermite polynomial of the degree $k-1$, see \cite{Suetin}, i.e.
\begin{eqnarray*}
	H_{k-1}(x)=(-1)^{k-1} e^{x^2}\cdot \dfrac{\dd^{k-1} }{\dd x^{k-1}}\,e^{-x^2},\quad x\in\R,\quad k\in\N. 
\end{eqnarray*}
For the general case,  $d\in\N$,  the sequence $(\lambda_{d, m})_{m\in\N}$ is an array of the numbers  
\begin{eqnarray*}
	\prod_{k = 1}^d \lambda_{k_j},\quad k_1,\ldots, k_d \in\N,
\end{eqnarray*}
that is indexed in the non-increasing order, and $(\psi_{d, m})_{m\in\N}$ is  correspondingly indexed sequence of the functions
\begin{eqnarray*}
	\prod_{k = 1}^d \psi_{k_j}(t_j),\quad k_1,\ldots, k_d \in\N.
\end{eqnarray*}
Observe that $\lambda_k$ and $\lambda_{d, m}$ belong to the interval $(0, 1)$ for all $k, d, m \in \N$. Moreover, it is not difficult to show that
$\sum_{k \in \N}\lambda_k = 1$ and $\sum_{m \in \N} \lambda_{d, m} = 1$ for every $d \in \N$.

For a given shape parameter $\gamma>0$ we consider the information complexity $n(d, \e)$ as a function depending on two variables
$d \in \N$, $\e \in (0, 1)$, and we are interested in the character of growth of $n(d, \e)$ for arbitrarily large $d$ and small $\e$. There are a lot of results for this and more general problems concerning tractability (see  \cite{FassHickWoz2} and \cite{SloanWoz}), where necessary and sufficient
conditions for the special upper bounds of $n(d, \e)$  were obtained in terms of input  parameters. Some of these results we will recall in the next section. 
There exist results about tractabilty within average case setting  (see \cite{ChenWang}, \cite{FassHickWoz1}, and \cite{Khart2}). This setting is also studied in \cite{KhartLim} with asymptotic approach  that more accurately detect behaviour of the quantity $n(d, \e)$ as $d\to\infty$. We are not aware any results with asymptotic approach for the worst case setting. In this paper, we do the first step in this direction. Namely, we will investigate asymptotics of $n(d, \e)$ as $d \to \infty$ for arbitrarily small fixed  $\varepsilon \in (0, 1)$ and for the case $\e=\e_d\to 0$ as $d\to\infty$. The latter setting seems to be new. 

We will use the following notation in the paper. The set of non-negative integers will be denoted by $\N_0$.   The indicator function $\mathds{1}(A)$ equals one if $A$ is true and it is zero if $A$ is false. For a finite set $B$ we denote by $\# B$ the number of elements in the set $B$. By 
\textit{distribution function} we mean a non-decreasing function $F$ defined on $\R$ that is right-continuous, $\lim\limits_{x \to +\infty}F(x) = 1$, and
$\lim\limits_{x \to -\infty}F(x) = 0$. We denote by $\Probab(C)$ the probability of an event $C$, and we donote by $\Expec Y$, and  $\Var Y$ the expectation and the variance of a random variable $Y$, respectively.

\section{Main results}

In this section we formulate  the main results of the paper. The proofs are provided in section~4, the necessary auxiliary tools are presented in section~3.

We consider the information complexity $n(d, \e)$ defined in the previous section. We first recall existing results for this quatity that is corresponded to the kernel \eqref{kernel-eqnarray} with any fixed $\gamma>0$. From \cite{FassHickWoz2} it is known (see also review in \cite{SloanWoz}) that \textit{quasi-polynomial tractability} (QPT) holds for $n(d,\e)$, i.e.  there  exist constants $C > 0$ and $\tau > 0$ such that
\begin{eqnarray}\label{QPT}
	n(d, \e) \leqslant C \exp\bigl\{\tau (1 + \ln d)(1 + \ln \e^{-1})\bigr\}, \quad\text{for all}\quad \e \in (0, 1), \,\, d \in \N.
\end{eqnarray}
Here the constant $\tau$  was determined as $\tau=2c_\omega$, where 
\begin{eqnarray}
	\label{const-omega-eqnarray}
	c_\omega \defeq \frac{1}{|\ln \omega|}.
\end{eqnarray}
The constant $c_\omega$ will be used here and throughout the paper.
  
Recall that there are two natural type of tractability that would be possible for $n(d,\e)$:
\begin{itemize}
	\item \textit{strong polynomial tractability} (SPT)  holds iff there exist constant $C > 0$ and $\tau > 0$ such that
	\begin{eqnarray*}
		n(d, \e) \leqslant C \e^{-\tau}, \quad\text{for all}\quad \e \in (0, 1), \,\, d \in \N;
	\end{eqnarray*}
	\item \textit{polynomial tractability} (PT) holds iff there exist constants $C > 0$, $\tau_1 > 0$ and $\tau_2 > 0$ such that
	\begin{eqnarray*}
		n(d, \e) \leqslant C \e^{-\tau_1} d^{\tau_2},\quad\text{for all}\quad \e \in (0, 1), \,\, d \in \N.
	\end{eqnarray*}
\end{itemize}
However, as it was shown in \cite{FassHickWoz2}, SPT and PT do not hold for the quantity $n(d,\e)$ for the considered case. Morever, \textit{exponential convergence and quasi-polynomial tractability} (EC-QPT) does not hold (see \cite{SloanWoz}), where EC-QPT, by its definition, holds iff there exist constants $C > 0$ and $\tau > 0$ such that
\begin{eqnarray*}
	n(d, \e) \leqslant C \exp\bigl\{\tau (1 + \ln d)(1 + \ln (1 + \ln \e^{-1}))\bigr\}, \quad \e \in (0, 1), \quad d \in \N.
\end{eqnarray*}  

We now formulate our asymptotical results that complement the estimate \eqref{QPT}. We first consider the case when $d\to\infty$ and $\e$ is fixed.

\begin{Theorem}\label{fixed_error_theorem}
    For any fixed $\e \in (0, 1)$
    \begin{eqnarray}
         \label{result1-eqnarray}
         n(d, \e) = \dfrac{d^{N_{\omega}(\varepsilon)}}{N_{\omega}(\varepsilon)!} + O\bigl(d^{N_{\omega}(\varepsilon)-1}\bigr),\quad d \to \infty,
    \end{eqnarray}
    where $N_{\omega}(\e)\defeq\max\bigl\{m\in\N_0: m< c_\omega |\ln \e^2|\bigr\}$.
\end{Theorem}

We see that the main factor $d^{N_{\omega}(\varepsilon)}$ in \eqref{result1-eqnarray} is approximately equal to $\exp\bigl\{2c_\omega\cdot \ln d\cdot |\ln \e|\bigr\}$. This is consistent with \eqref{QPT} (with $\tau=2c_\omega$) for large $d$ and small $\e$.

We turn to the case when $d\to\infty$ and $\e=\e_d\to 0$. Without loss of generality concerning the rate of convergence of $(\e_d)_{d\in\N}$, we  obtain the logarithmic asymptotics of $n(d,\e_d)$ as $d\to\infty$.

\begin{Theorem}\label{logarithmic_nonfixed_error_theorem}
    Let $(\e_d)_{d\in\N}$  be a sequence from $(0,1)$ such that $\e_d\to 0$, $d\to\infty$.  Then
    \begin{eqnarray}
        \label{result2-eqnarray}
        \ln n(d, \varepsilon_d) = d\cdot\ln\Bigl(1+\tfrac{c_\omega|\ln \e_d^2|}{d}\Bigr) +
        c_\omega|\ln \e_d^2|\cdot\ln\Bigl(1+\tfrac{d}{c_\omega|\ln \e_d^2|}\Bigr)+R_d, \quad d \to \infty.
    \end{eqnarray}
where $(R_d)_{d\in\N}$ is a sequence such that
\begin{eqnarray*}
	\dfrac{R_d}{|\ln \e_d^2|\cdot\ln\Bigl(1+\tfrac{d}{c_\omega|\ln \e_d^2|}\Bigr)}\to0,\quad d\to\infty.
\end{eqnarray*}
\end{Theorem}

For large $d$ the asymptotics \eqref{result2-eqnarray}  is  consistent with \eqref{QPT} too. Indeed, if $d / |\ln \e_d^2|\to \infty$ as $d\to\infty$, then for large $d$ the quantity $\ln n(d, \varepsilon_d)$ is approximately equal to 
\begin{eqnarray*}
	c_\omega |\ln \e_d^2|+c_\omega |\ln \e_d^2| \ln d.
\end{eqnarray*} 
Here the second term is the main as in \eqref{QPT}. If $ |\ln \e_d^2|/d\to \infty$ as $d\to\infty$, then  $\ln n(d, \varepsilon_d)$ is approximately equal to
\begin{eqnarray*}
	d \ln\bigl(c_\omega|\ln \e_d^2|\bigr)+  d
\end{eqnarray*} 
for sufficiently large $d$. So the first term is the main and for it we have
\begin{eqnarray*}
	\dfrac{d \ln\bigl(c_\omega|\ln \e_d^2|\bigr)}{ c_\omega |\ln \e_d^2|\ln d}=  \dfrac{d}{ c_\omega |\ln \e_d^2|}\cdot \dfrac{\ln\Bigl(\tfrac{c_\omega|\ln \e_d^2|}{d}\Bigr)+\ln d}{ \ln d}=\dfrac{\ln\Bigl(\tfrac{c_\omega|\ln \e_d^2|}{d}\Bigr)}{\tfrac{c_\omega|\ln \e_d^2|}{d}}\cdot \dfrac{1}{\ln d}+ \dfrac{d}{ c_\omega |\ln \e_d^2|}.
\end{eqnarray*}
It is seen that, by the assumption,  this sequence goes to zero as $d\to\infty$. In particular, $d \ln\bigl(c_\omega|\ln \e_d^2|\bigr)$  is  bounded by $c_\omega |\ln \e_d^2|\ln d$ for all sufficiently large $d$. Hence the same estimate holds for  $\ln n(d, \varepsilon_d)$ similarly to \eqref{QPT} with $\tau=2c_\omega$.

\section{Auxiliary results}

In this section, we formulate and prove some statements that are needed to obtain the main results. 

\begin{Lemma} \label{lem_first_step}
    Let $(U_j)_{j \in \N}$ be a sequence of independent random variables with the following distribution:
    \begin{eqnarray}\label{lem_first_step_Uj}
        \Probab\bigl(U_j = k|\ln \omega|\bigr) = (1 - \omega)\,\omega^k, \quad k \in \N_0, \quad j\in\N.
    \end{eqnarray}
Then for any $d\in\N$ and $\e\in (0,1)$ the quantity $n(d, \e)$ admits the following representation
    \begin{eqnarray*} 
        n(d, \e) =e^{d|\ln (1 - \omega)|}\,\Expec \biggl[\exp\Bigl\{\sum_{j = 1}^d  U_j\Bigr\}\id\Bigl(\sum_{j = 1}^d U_j< |\ln \varepsilon^2|\Bigr)\biggr].
    \end{eqnarray*}
\end{Lemma}
\textbf{Proof of Lemma \ref{lem_first_step}.}\quad Let us fix any $d\in\N$ and $\e\in(0,1)$.  Using  \eqref{approx-compl-repr1-eqnarray}, we have
\begin{eqnarray*}
    n(d, \varepsilon)=\min\bigl\{n \in \N: \lambda_{d, n + 1} \leqslant \e^2 \lambda_{d, 1}\bigr\}
    = \#\{n\in\N: \lambda_{d, n + 1} > \e^2 \lambda_{d, 1}\}.
\end{eqnarray*}
According to the known product structure of  $\lambda_{d, m}$, $m\in\N$, and due to the sharp expressions \eqref{onedim-eigenvalue-eqnarray} for $\lambda_k$, $k\in\N$, we have
\begin{eqnarray*}
	n(d, \varepsilon)
	&=&\#\Big\{(k_1, \ldots, k_d) \in \N^d: \prod_{j= 1}^d \lambda_{k_j} > \e^2 \lambda_{d, 1}\Big\} \\
	&=&\#\Big\{(k_1, \ldots, k_d) \in \N^d: \prod_{j = 1}^d \bigl((1 - \omega)\,\omega^{k_j - 1}\bigr) > \e^2 (1 - \omega)^d\Big\} \\
	&=&\#\Big\{(k_1, \ldots, k_d) \in \N^d_0: \prod_{j = 1}^d \omega^{k_j } > \e^2\Big\} \\
	&=&\#\Big\{(k_1, \ldots, k_d) \in \N_0^d: \sum_{j = 1}^d k_j |\ln \omega| < |\ln \e^2|\Big\}.
\end{eqnarray*}
We write
\begin{eqnarray*}
    n(d, \varepsilon)&=&
    \sum_{(k_1, \ldots, k_d) \in \N_0^d}\id\Big(\sum_{j = 1}^d k_j |\ln \omega| < |\ln \e^2|\Big) \\ &=&e^{d|\ln (1 - \omega)|}
    \sum_{(k_1, \ldots, k_d) \in \N_0^d} \biggl[\exp\Bigl\{\sum_{j = 1}^d k_j |\ln \omega|\Bigr\}\prod_{j = 1}^d \bigl((1 - \omega)\,\omega^{k_j}\bigr)\id\Big(\sum_{j = 1}^d k_j |\ln \omega| < |\ln \e^2|\Bigr)\biggr].
\end{eqnarray*}
By definition \eqref{lem_first_step_Uj}, we obtain
\begin{eqnarray*}
    n(d, \varepsilon)
    =e^{d|\ln (1 - \omega)|}\Expec \biggl[\exp\Big\{\sum_{j = 1}^d U_j\Big\} \id\Big(\sum_{j = 1}^d U_j < |\ln \e^2|\Big)\biggr],
\end{eqnarray*}
as required. \quad $\Box$\\

Let us define the moment generating function for random variables $U_j$, $j\in\N$:
\[
M(\nu) \defeq \Expec \exp\{\nu U_j\}= \sum_{k=0}^{\infty} e^{\nu|\ln \omega|k}(1 - \omega)\,\omega^k=\sum_{k=0}^{\infty}(1 - \omega)\, \omega^{(1-\nu)k},\quad \nu\in\R,
\]
The next lemma states some properties of this function and its derivatives.

\begin{Lemma} \label{lem_second_step} $1)$ For any $\nu<1$  the quantities $M(\nu)$, $M'(\nu)$, and $M''(\nu)$ are finite and admit representations
        \begin{eqnarray*}
            M(\nu) = \frac{1 - \omega}{1 - \omega^{1 - \nu}} ,\qquad 
            M'(\nu) = \frac{(1 - \omega)|\ln \omega|\,\omega^{1 - \nu}}{(1 - \omega^{1 -\nu})^2},\qquad   M''(\nu) =    (1 - \omega)|\ln \omega|^2\cdot\dfrac{(1 + \omega^{1 - \nu})\,\omega^{1 - \nu}}{(1 - \omega^{1 - \nu})^3}.
        \end{eqnarray*}
         $2)$ For any $\alpha > 0$ the equation $\alpha = M'(\nu)/M(\nu)$ has the solution $\nu_{\alpha}<1$ expressed by the formula
        \begin{eqnarray}\label{eq_nu_alpha}
            \nu_\alpha =  1- c_\omega \ln\bigl(1+\tfrac{1}{c_\omega\alpha}\bigr),
        \end{eqnarray}
and then
\begin{eqnarray}\label{eq_M}
	M(\nu_\alpha) = (1 - \omega)(1+c_\omega \alpha),\qquad \dfrac{M'(\nu_\alpha)}{M(\nu_{\alpha})}=\alpha,\qquad \dfrac{M''(\nu_\alpha)}{M(\nu_{\alpha})}- \biggl(\dfrac{M'(\nu_\alpha)}{M(\nu_{\alpha})}\biggr)^2= |\ln \omega|\cdot (1+c_\omega \alpha)\cdot \alpha.
\end{eqnarray}

\end{Lemma}
\textbf{Proof of Lemma \ref{lem_second_step}.}\quad  $1)$ Since $\omega \in (0, 1)$, we get
\begin{eqnarray*}
   M(\nu) = (1 - \omega)\sum_{k=0}^{\infty} \omega^{(1-\nu)k}
    = \dfrac{1 - \omega}{1 - \omega^{1 - \nu}}<\infty, \quad \nu<1.
\end{eqnarray*}
This yields
\begin{eqnarray*}
    M'(\nu) =\biggl(\dfrac{1 - \omega}{1 - \omega^{1 - \nu}}\biggr)'=  
    \dfrac{ (1 - \omega)|\ln \omega|\,\omega^{1 - \nu}}{(1 - \omega^{1 - \nu})^2}<\infty,\quad \nu<1.
\end{eqnarray*}
Next, we have
\begin{eqnarray*}
	M''(\nu) &=& \biggl(\dfrac{ (1 - \omega)|\ln \omega|\,\omega^{1 - \nu}}{(1 - \omega^{1 - \nu})^2}\biggr)' \\
	&=& (1 - \omega)|\ln\omega|\cdot\dfrac{|\ln \omega|\,\omega^{1 - \nu}\cdot(1 - \omega^{1 - \nu})^2 + \omega^{1 - \nu}\cdot 2(1 - \omega^{1 - \nu})\,\omega^{1 - \nu}\,|\ln \omega|\,  }
	{(1 - \omega^{1 - \nu})^4} \\
	&=&(1 - \omega)|\ln\omega|\cdot\dfrac{|\ln \omega|\,\omega^{1 - \nu}(1 - \omega^{1 - \nu}) + 2|\ln \omega|\, \omega^{2(1 - \nu)}}
	{(1 - \omega^{1 - \nu})^3} \\ &=&
	(1 - \omega)|\ln \omega|^2\cdot\dfrac{(1 + \omega^{1 - \nu})\,\omega^{1 - \nu}}{(1 - \omega^{1 - \nu})^3}<\infty,\quad \nu<1.
\end{eqnarray*}

$2)$  We fix any $\alpha>0$ and consider the equation $\alpha = M'(\nu)/M(\nu)$ with $\nu<1$. On account of the obtained formulas, we have
\begin{eqnarray*}
	\dfrac{M'(\nu)}{M(\nu)} = \frac{|\ln \omega|\cdot \omega^{1 - \nu}}{1 - \omega^{1 - \nu}},\quad \nu<1.
\end{eqnarray*}
So the equation gets the form
\begin{eqnarray*}
	\alpha =\dfrac{|\ln \omega|\cdot\omega^{1 - \nu}}{1 - \omega^{1 - \nu}},
\end{eqnarray*}
that is equivalent to
\begin{eqnarray*}
    \omega^{1 - \nu} = \frac{\alpha}{ |\ln \omega|+\alpha}\quad \text{and}\quad (1 - \nu)|\ln \omega| = \ln\frac{ |\ln \omega|+\alpha}{\alpha}.
\end{eqnarray*}
The last logarithm equals $\ln\bigl(1+ \tfrac{1}{c_\omega \alpha}\bigr)$. So we easily get the unique solution \eqref{eq_nu_alpha} that is less than $1$. Thus the second equality in \eqref{eq_M} is obvious. We now check the first:
\begin{eqnarray*}
	M(\nu_{\alpha})= \dfrac{1 - \omega}{1 - \omega^{1 - \nu_{\alpha}}}= \dfrac{1 - \omega}{1 - \tfrac{\alpha}{|\ln \omega|+\alpha}}= \dfrac{1 - \omega}{ \tfrac{|\ln \omega|}{|\ln \omega|+\alpha}}=(1-\omega) \Bigl(1+\tfrac{\alpha}{|\ln \omega|}\Bigr)=(1 - \omega)(1+c_\omega \alpha).
\end{eqnarray*}
Next, we prove the third equality in \eqref{eq_M}. By the above, we have
\begin{eqnarray*}
	\dfrac{M''(\nu_\alpha)}{M(\nu_{\alpha})}&=&\dfrac{|\ln \omega|\,(1 + \omega^{1 - \nu_\alpha})}{1 - \omega^{1 - \nu_\alpha}}\cdot\dfrac{M'(\nu_\alpha)}{M(\nu_{\alpha})}\\
	&=&|\ln \omega|\cdot \Bigl(1+ \frac{\alpha}{|\ln \omega|+\alpha}\Bigr)\cdot(1+c_\omega \alpha)\cdot \alpha \\
	&=& \Bigl(|\ln \omega|+ \frac{\alpha}{1+c_\omega\alpha}\Bigr)\cdot (1+c_\omega \alpha)\cdot\alpha\\
	&=&|\ln \omega|\cdot (1+c_\omega \alpha)\cdot\alpha +\alpha^2.
\end{eqnarray*}
Thus we obtain
\begin{eqnarray*}
	\dfrac{M''(\nu_\alpha)}{M(\nu_{\alpha})}- \biggl(\dfrac{M'(\nu_\alpha)}{M(\nu_{\alpha})}\biggr)^2=\dfrac{M''(\nu_\alpha)}{M(\nu_{\alpha})}-\alpha^2= |\ln \omega|\cdot (1+c_\omega \alpha)\cdot \alpha.\quad \Box
\end{eqnarray*}

\begin{Lemma}\label{lem_third_step}
	 Let for every $\nu<1$ $(U_j^{(\nu)})_{j \in \N}$ be a sequence of independent random variables with the following distribution:
    \begin{eqnarray*}
        \Probab\bigl(U^{(\nu)}_j = k|\ln \omega|\bigr) = \dfrac{(1 - \omega)\,\omega^{(1 - \nu)k}}{M(\nu)}, \quad k \in \N_0, \quad j\in\N.
    \end{eqnarray*}
 Then for any $d\in\N$,  $\e\in (0,1)$, and $\nu<1$ the quantity $n(d, \e)$ admits the following representation
\begin{eqnarray} \label{eq_third_step}
	n(d, \e) =e^{d|\ln (1 - \omega)|} M(\nu)^d\,\Expec \biggl[\exp\Big\{(1-\nu)\sum_{j = 1}^d U_j^{(\nu)}\Big\} \id\Big(\sum_{j = 1}^d U_j^{(\nu)} < |\ln \e^2|\Big)\biggr].
\end{eqnarray}
\end{Lemma}
\textbf{Proof of Lemma \ref{lem_third_step}.}\quad Let us fix arbitrary $d\in\N$,  $\e\in (0,1)$, and $\nu<1$. Observe that
\begin{eqnarray*}
    \Probab\bigl(U^{(\nu)}_j = k|\ln \omega|\bigr) = \dfrac{e^{k|\ln \omega|\nu}}{M(\nu)}\, (1 - \omega)\,\omega^{k} = \dfrac{e^{k|\ln \omega|\nu}}{M(\nu)}\,\Probab\bigl(U_j = k|\ln \omega|\bigr),\quad k\in\N_0,\quad j\in\N.
\end{eqnarray*}
Hence, due to the independence, for any $k_1,\ldots, k_d\in\N_0$ we have 
\begin{eqnarray*}
	\Probab\bigl(U^{(\nu)}_1 = k_1|\ln \omega|,\ldots, U^{(\nu)}_d = k_d|\ln \omega| \bigr)  &=& 
	\Probab\bigl(U^{(\nu)}_1 = k_1|\ln \omega|\bigr)\cdot\ldots \cdot\Probab\bigl( U^{(\nu)}_d = k_d|\ln \omega| \bigr)\\
	&=&\dfrac{e^{k_1|\ln \omega|\nu}}{M(\nu)}\,\Probab\bigl(U_1 = k_1|\ln \omega|\bigr)\cdot\ldots \cdot\dfrac{e^{k_d|\ln \omega|\nu}}{M(\nu)}\,\Probab\bigl( U_d = k_d|\ln \omega| \bigr)\\
	&=&\dfrac{e^{(k_1+\ldots+k_d)|\ln \omega|\nu}}{M(\nu)^d}\,  \Probab\bigl(U_1 = k_1|\ln \omega|,\ldots, U_d = k_d|\ln \omega| \bigr).
\end{eqnarray*}

Let $F_d$ and $F_d^{(\nu)}$ be the distribution functions of the sums $\sum_{j =1}^d U_j$ and $\sum_{j =1}^d U_j^{(\nu)}$, respectively.
Then for any $x\geqslant 0$ we get
\begin{eqnarray*}
	F_d^{(\nu)}(x)&=&\Probab\Bigl(\sum_{j=1}^{d} U_j^{(\nu)}\leqslant x\Bigr)\\&=&\sum_{\substack{k_1,\ldots, k_d\in\N_0:\\ k_1+\ldots+k_d\leqslant x}}  \dfrac{e^{(k_1+\ldots+k_d)|\ln \omega|\nu}}{M(\nu)^d}\,  \Probab\bigl(U_1 = k_1|\ln \omega|,\ldots, U_d = k_d|\ln \omega| \bigr)\\
	&=& \int_{0}^{x}\dfrac{e^{\nu u}}{M(\nu)^d}\, \dd F_d(u).
\end{eqnarray*}
Applying  Lemma \ref{lem_first_step} together with this formula,  we obtain the needed representation:
\begin{eqnarray*}
    n(d, \varepsilon)&=&e^{d|\ln (1 - \omega)|}\,\Expec \biggl[\exp\Bigl\{\sum_{j = 1}^d  U_j\Bigr\}\id\Bigl(\sum_{j = 1}^d U_j < |\ln \varepsilon^2|\Bigr)\biggr]\\
    &=& e^{d|\ln (1 - \omega)|}\int_{\R} e^x \id\bigl(x < |\ln \varepsilon^2|\bigr)\,\dd F_d(x)\\
    &=& e^{d|\ln (1 - \omega)|}M(\nu)^d\int_{\R} e^{(1-\nu)x}\id\bigl(x < |\ln \varepsilon^2|\bigr)\, \dd F_d^{(\nu)}(x)\\
    &=&e^{d|\ln (1 - \omega)|} M(\nu)^d\,\Expec \biggl[\exp\Big\{(1-\nu)\sum_{j = 1}^d U_j^{(\nu)}\Big\} \id\Big(\sum_{j = 1}^d U_j^{(\nu)} < |\ln \e^2|\Big)\biggr].\quad \Box
\end{eqnarray*}

Note that formula \eqref{eq_third_step} can be written in the following form:
\begin{eqnarray*}
	n(d, \e) =e^{d|\ln (1 - \omega)|} M(\nu)^d\,e^{(1-\nu)|\ln\e^2|}\Expec \biggl[\exp\Big\{(1-\nu)\Bigl(\sum_{j = 1}^d U_j^{(\nu)}-|\ln\e^2|\Bigr)\Bigr\} \id\Big(\sum_{j = 1}^d U_j^{(\nu)} < |\ln \e^2|\Big)\biggr].
\end{eqnarray*}
It is easily seen that the expectation less than $1$ due to $\nu<1$. Hence the part before this dominates the growth of  the quantity $n(d, \e)$. Let us choose  $\nu<1$ to  minimize this part. It is equivalent to the minimization of the function $\nu\mapsto M(\nu)^d\,e^{-\nu|\ln\e^2|}$ and, consequently, of the function $\nu\mapsto\ln M(\nu)-\nu\tfrac{|\ln\e^2|}{d}$. It is easily seen that for any $d\in\N$ and $\e\in(0,1)$ the latter function has the minimum at the point $\nu_{d, \e}$ such that 
\begin{eqnarray*}
	\dfrac{M'(\nu_{d, \e})}{M(\nu_{d, \e})}=\dfrac{|\ln \e^2|}{d}.
\end{eqnarray*}
Applying Lemma \ref{lem_second_step}  with $\alpha=\tfrac{|\ln \e^2|}{d}$ (and $\nu_{\alpha}=\nu_{d, \e}$), we find
 \begin{eqnarray}
      \label{eq_nu_d_eps}
	\nu_{d, \e} = 1 - c_\omega\ln\Bigl(1+\tfrac{d}{c_\omega|\ln \e^2|}\Bigr),\quad d\in\N,\,\, \e\in(0,1). 
\end{eqnarray}
Substituting this in \eqref{eq_third_step} and using the first equality of \eqref{eq_M}, we have
\begin{eqnarray}\label{conc_Mnude1minusomega}
	 e^{d|\ln (1 - \omega)|} M(\nu_{d,\e})^d=  e^{d|\ln (1 - \omega)|}(1 - \omega)^d\Bigl(1+\tfrac{c_\omega|\ln \e^2|}{d} \Bigr)^d=\Bigl(1+\tfrac{c_\omega|\ln \e^2|}{d} \Bigr)^d,
\end{eqnarray}
and thus
\begin{eqnarray}\label{eq_fourth_step}
	n(d, \e) = \Bigl(1+\tfrac{c_\omega|\ln \e^2|}{d} \Bigr)^d\,
	\Expec \biggl[\exp\Big\{c_\omega\ln\Bigl(1+\tfrac{d}{c_\omega|\ln \e^2|}\Bigr)S_{d, \e}\Big\}
	\id\Big(S_{d, \e} < |\ln \e^2|\Big)\biggr],\quad d\in\N,\,\, \e\in(0,1),
\end{eqnarray}
where we define
\begin{eqnarray*}
	S_{d,\e}\defeq\sum_{j = 1}^d U_j^{(\nu_{d,\e})},\quad d\in\N,\,\, \e\in(0,1). 
\end{eqnarray*}

Now we consider the asymptotic behaviour of distributions of $S_{d, \e}$ as $d\to\infty$ for the case of fixed $\e \in (0, 1)$ and for the case, when $\e=\e_d\to 0$ as $d\to\infty$.

\begin{Proposition}\label{prop-poisson}
    For any fixed $m\in\N_0$ and $\e \in (0, 1)$ 
    \begin{eqnarray*}
        \Probab\bigl(S_{d, \e} = m|\ln \omega|\bigr)=
        \exp\bigl\{-c_\omega|\ln \e^2|\bigr\}\cdot\dfrac{(c_\omega|\ln \e^2|)^m}{m!}+O\bigl(\tfrac{1}{d}\bigr),\quad d\to\infty.
    \end{eqnarray*}
\end{Proposition}
\textbf{Proof of Proposition \ref{prop-poisson}.}\quad 
Let us fix $m\in\N_0$ and $\e \in (0, 1)$. Since $S_{d, \e}$ is a sum of $d$ independent and identically distributed random variables, we have
\begin{eqnarray*}
    \Probab\bigl(S_{d, \e} = m|\ln \omega|\bigr) &=&
    \sum_{\substack{k_1, \ldots, k_d \in \N_0:\\ k_1 + \ldots + k_d = m}}\Probab\Big(U_1^{(\nu_{d, \e})} = k_1|\ln \omega|,\ldots,U_d^{(\nu_{d, \e})} = k_d|\ln \omega|\Big)\\
    &=&
    \sum_{\substack{k_1, \ldots, k_d \in \N_0: \\ k_1 + \ldots + k_d = m}}\Probab\Big(U_1^{(\nu_{d, \e})} = k_1|\ln \omega|\Big)\cdot\ldots\cdot\Probab
    \Big(U_d^{(\nu_{d, \e})} = k_d|\ln \omega|\Big) \\ &=&
    \sum_{\substack{k_1, \ldots, k_d \in \N_0:\\ k_1 + \ldots + k_d = m}}\prod_{j = 1}^d \dfrac{(1 - \omega)\,\omega^{k_j(1 - \nu_{d, \e})}}{M(\nu_{d, \e})} \\ &=&
    \sum_{\substack{k_1, \ldots, k_d \in \N_0:\\ k_1 + \ldots + k_d = m}}\frac{(1 - \omega)^d\omega^{m(1 - \nu_{d, \e})}}{M(\nu_{d,\e})^d}.
\end{eqnarray*}
It is known that the equation $k_1 + \ldots  + k_d = m$ with respect to non-negative integers $k_1, \ldots, k_d$ has 
$\binom{m + d - 1}{m}$ solutions (see, \cite{Feller}, p. 38). Therefore
\begin{eqnarray*}
	\Probab\bigl(S_{d, \e} = m|\ln \omega|\bigr) =\binom{m + d - 1}{m}\frac{(1 - \omega)^d\,\omega^{m(1 - \nu_{d, \e})}}{M(\nu_{d,\e})^d},\quad d\in\N.
\end{eqnarray*}
We represent
\begin{eqnarray*}
	\binom{m + d - 1}{m}= \frac{d \cdot \ldots \cdot (d + m - 1)}{m!}= \frac{1 \cdot \ldots \cdot \left(1 + \tfrac{m - 1}{d}\right)}{m!}\cdot d^m.
\end{eqnarray*}
Note that  the product in the fraction is assumed to be $1$ for the case $m=0$. Next, on account of \eqref{eq_nu_d_eps} and also \eqref{const-omega-eqnarray}, we have
\begin{eqnarray}\label{conc_omegam1nu}
	\omega^{m(1 - \nu_{d, \e})}=e^{ m (1 - \nu_{d, \e})\ln \omega}=\Bigl(1+ \tfrac{d}{c_\omega|\ln \e^2|}\Bigr)^{m c_\omega\ln\omega}=\Bigl(1+ \tfrac{d}{c_\omega|\ln \e^2|}\Bigr)^{-m }=\Bigl(\tfrac{c_\omega|\ln \e^2|}{d}\Bigr)^m \Bigl(1+ \tfrac{c_\omega|\ln \e^2|}{d}\Bigr)^{-m}.
\end{eqnarray}
Due to \eqref{conc_Mnude1minusomega}, we get
\begin{eqnarray}\label{conc_1omegaM}
	\dfrac{(1 - \omega)^d}{M(\nu_{d,\e})^d}=e^{-d|\ln (1 - \omega)|} M(\nu_{d,\e})^{-d}=\Bigl(1+\tfrac{c_\omega|\ln \e^2|}{d} \Bigr)^{-d}.
\end{eqnarray}
Thus
\begin{eqnarray*}
	\Probab\bigl(S_{d, \e} = m|\ln \omega|\bigr) =\dfrac{1 \cdot \ldots \cdot \left(1 + \tfrac{m - 1}{d}\right)}{m!} \bigl(c_\omega|\ln \e^2|\bigr)^m \Bigl(1+ \tfrac{c_\omega|\ln \e^2|}{d}\Bigr)^{-m} \Bigl(1+\tfrac{c_\omega|\ln \e^2|}{d} \Bigr)^{-d},\quad d\in\N.
\end{eqnarray*}
It is easily seen that for fixed $m\in\N_0$ and $\e\in(0,1)$
\begin{eqnarray*}
	1 \cdot \ldots \cdot \left(1 + \tfrac{m - 1}{d}\right)=1+O\bigl(\tfrac{1}{d}\bigr),\qquad \Bigl(1+ \tfrac{c_\omega|\ln \e^2|}{d}\Bigr)^{-m}=1+O\bigl(\tfrac{1}{d}\bigr),\quad d\to\infty.
\end{eqnarray*}
Next, observe that for $d\to\infty$
\begin{eqnarray*}
	\Bigl(1+\tfrac{c_\omega|\ln \e^2|}{d} \Bigr)^{-d}= \exp\Bigl\{-d \ln\Bigl(1 + \tfrac{c_\omega|\ln \e^2|}{d}\Bigr)\Bigr\}=\exp\bigl\{-c_\omega|\ln \e^2|+ O\bigl(\tfrac{1}{d}\bigr) \bigr\}.
\end{eqnarray*}
So we have
\begin{eqnarray}\label{conc_1plusedd}
		\Bigl(1+\tfrac{c_\omega|\ln \e^2|}{d} \Bigr)^{-d}=\exp\bigl\{-c_\omega|\ln \e^2| \bigr\}\Bigl(1+ O\bigl(\tfrac{1}{d}\bigr)\Bigr),\quad d\to\infty.
\end{eqnarray}

Therefore
\begin{eqnarray*}
	\Probab\bigl(S_{d, \e} = m|\ln \omega|\bigr) &=&\dfrac{1+O\bigl(\tfrac{1}{d}\bigr)}{m!} \bigl(c_\omega|\ln \e^2|\bigr)^m \Bigl(1+ O\bigl(\tfrac{1}{d}\bigr)\Bigr)\cdot \exp\bigl\{-c_\omega|\ln \e^2| \bigr\}\Bigl(1+ O\bigl(\tfrac{1}{d}\bigr)\Bigr)\\
	&=& \exp\bigl\{-c_\omega|\ln \e^2|\bigr\}\cdot\dfrac{(c_\omega|\ln \e^2|)^m}{m!} \Bigl(1+ O\bigl(\tfrac{1}{d}\bigr)\Bigr)\\
	&=& \exp\bigl\{-c_\omega|\ln \e^2|\bigr\}\cdot\dfrac{(c_\omega|\ln \e^2|)^m}{m!} + O\bigl(\tfrac{1}{d}\bigr),\quad d\to\infty.\quad \Box\\
\end{eqnarray*}

\begin{Proposition}\label{normal_proposition}
	Let $(\e_d)_{d\in\N}$  be a sequence from $(0,1)$ such that $\e_d\to 0$, $d\to\infty$.  Then
	\begin{eqnarray}\label{normal_proposition_conv}
		\sup_{x \in \R}\Biggl|\Probab\Biggl(\dfrac{S_{d, \e_d} - \Expec S_{d, \e_d}}{\sqrt{\Var S_{d, \e_d}}}\leqslant x\Biggr) - \Phi(x)\Biggr| \to 0, \quad d \to \infty,
	\end{eqnarray}
	where $\Phi$ is the distribution function of standard normal law:
	\begin{eqnarray*}
		\Phi(x)\defeq \dfrac{1}{\sqrt{\pi}}\int_{-\infty}^x e^{-u^2/2} \dd u,\quad x\in\R.
	\end{eqnarray*}
\end{Proposition}
\textbf{Proof of Proposition \ref{normal_proposition}.}\quad Suppose that the sequence  $(\e_d)_{d\in\N}$ is given. We define 
\begin{eqnarray}\label{def_adbd}
	a_d\defeq \Expec S_{d, \e_d},\quad\text{and} \quad b_d\defeq \sqrt{\Var S_{d, \e_d} },\quad d\in\N.
\end{eqnarray}
We will prove the convergence \eqref{normal_proposition_conv} using the method of characteristic functions (see \cite{Lukacs}). We introduce characteristic functions for the fractions in \eqref{normal_proposition_conv}:
\begin{eqnarray*}
	H_d(t) = \Expec \exp\Bigl\{  it \tfrac{S_{d, \e_d}-a_d}{b_d}\Bigr\},\quad t\in\R,\quad d\in\N.
\end{eqnarray*}
Since $S_{d, \e_d}$ is the sum of $d$ independent and identically distributed random variables $U_{j}^{(\nu_{d, \e_d})}$, we have 
\begin{eqnarray*}
	H_d(t) = h_d\bigl(\tfrac{t}{b_d}\bigr)^d\exp\bigl\{-it\tfrac{a_d}{b_d}\bigr\},\quad t\in\R,\quad d\in\N,
\end{eqnarray*}
where $h_d$ is the characteristic function of $U_{j}^{(\nu_{d, \e_d})}$, $d\in\N$, i.e.
\begin{eqnarray*}
	h_d(t)\defeq \Expec \exp\bigl\{  it U_{j}^{(\nu_{d, \e_d})}\bigr\},\quad t\in\R,\quad d\in\N.
\end{eqnarray*}
Let us find the convenient expression of $h_d(t)$ for any $t\in\R$ and $d\in\N$:
\begin{eqnarray*}
     h_d(t) = \sum_{k = 0}^{\infty} \dfrac{(1 - \omega)\omega^{k(1 - \nu_{d, \e_d})}}{M(\nu_{d, \e_d})}\, e^{itk|\ln \omega|}=\dfrac{1 - \omega}{M(\nu_{d, \e_d})}\sum_{k = 0}^{\infty} \bigl(\omega^{1 - \nu_{d, \e_d}}\cdot e^{it|\ln \omega|}\bigr)^k = \dfrac{\tfrac{1 - \omega}{M(\nu_{d, \e_d})}}{1 - \omega^{1 - \nu_{d, \e_d}}\cdot e^{it|\ln \omega|}}.
\end{eqnarray*}
From \eqref{conc_omegam1nu} (with $m=1$)  and \eqref{conc_1omegaM} (with $d=1$ in the powers) we obtain
\begin{eqnarray*}
	h_d(t)= \dfrac{\Bigl(1+\tfrac{c_\omega|\ln \e^2_d|}{d}\Bigr)^{-1}}{1-  \Bigl(1+ \tfrac{d}{c_\omega|\ln \e^2_d|}\Bigr)^{-1} e^{it|\ln \omega|}}= \dfrac{\tfrac{d}{d+c_\omega|\ln \e^2_d|}}{1-   \tfrac{c_\omega|\ln \e^2_d|}{d+c_\omega|\ln \e^2_d|}\, e^{it|\ln \omega|}}=\Bigl(1-  \tfrac{c_\omega|\ln \e^2_d|}{d}\, \bigl(e^{it|\ln \omega|}-1\bigr)\Bigr)^{-1}
\end{eqnarray*}
for any $t\in\R$ and $d\in\N$. We now return to $H_d$:
\begin{eqnarray*}
    H_d(t) = \Bigl(1-  \tfrac{c_\omega|\ln \e^2_d|}{d}\, \bigl(   \exp\bigl\{\tfrac{it}{b_d}|\ln \omega|\bigr\}-1\bigr)\Bigr)^{-d} \exp\bigl\{-it\tfrac{a_d}{b_d}\bigr\},\quad t\in\R,\quad d\in\N.
\end{eqnarray*}
Let us consider the sequences $a_d$ and $b_d$, $d\in\N$. For $a_d$ we have
\begin{eqnarray*}
		a_d=\Expec S_{d, \e_d} =\Expec\Bigl(\sum_{j = 1}^d U_j^{(\nu_{d,\e_d})}\Bigr)=d\,\Expec U_1^{(\nu_{d,\e_d})},\quad d\in\N,
\end{eqnarray*}
due to the identical distributions of $U_j^{(\nu_{d,\e_d})}$. For $b_d$, due to the independence of $U_j^{(\nu_{d,\e_d})}$, we get
\begin{eqnarray*}
	 b_d^2=  \Var S_{d, \e_d}=\Var\Bigl(\sum_{j = 1}^d U_j^{(\nu_{d,\e_d})}\Bigr)=d\,\Var  U_1^{(\nu_{d,\e_d})},\quad d\in\N.
\end{eqnarray*}
Recall that, by definition, the distribution of $U_j^{(\nu_{d,\e_d})}$ can be considered as the Laplace transform of the distribution of $U_j$. From the theory of large deviations, it is known (see \cite{Borovkov}, p. 243) that  
\begin{eqnarray*}
   \Expec U_1^{(\nu_{d,\e_d})} =
       \frac{M'(\nu_{d, \e_d})}{M(\nu_{d, \e_d})}, \quad\text{and}\quad
   \Var U_1^{(\nu_{d,\e_d})} =
       \frac{M''(\nu_{d, \e_d})}{M(\nu_{d,\e_d})}-\biggl(\frac{M'(\nu_{d, \e_d})}{M(\nu_{d, \e_d})}\biggr)^2,\quad d\in\N.
\end{eqnarray*}
These equalities can be also checked directly using Lemma \ref{lem_second_step} and elementary notes about series (see \cite{KhartLim}, Lemma 1). Next, according to \eqref{eq_M} with $\alpha = \tfrac{|\ln \e_d^2|}{d}$ and $\nu_\alpha = \nu_{d, \e}$, we obtain
\begin{eqnarray*}
	\Expec U_1^{(\nu_{d,\e_d})} = \tfrac{|\ln \e_d^2|}{d},\quad \text{and} \quad \Var U_1^{(\nu_{d,\e_d})}= |\ln \omega|\cdot\Bigl(1+ \tfrac{c_\omega|\ln \e_d^2|}{d}\Bigr)\cdot \tfrac{|\ln \e_d^2|}{d},\quad d\in\N.
\end{eqnarray*}
Thus we have
\begin{eqnarray}\label{conc_adbd}
	a_d= |\ln \e_d^2|,\quad \text{and}\quad b_d^2= |\ln \omega|\cdot\Bigl(1+ \tfrac{c_\omega|\ln \e_d^2|}{d}\Bigr)\cdot |\ln \e_d^2|,\quad d\in\N.
\end{eqnarray}
Since $\e_d\to 0$ as $d\to\infty$, it is seen that $a_d\to \infty$ and $b_d\to\infty$ as $d\to\infty$. We also conclude that
\begin{eqnarray*}
    \label{eq_ad_div_dbd}
    \dfrac{a_d}{db_d} &=&
     \dfrac{|\ln \e_d^2|}{d\sqrt{|\ln \omega|\cdot\Bigl(1+ \tfrac{c_\omega|\ln \e_d^2|}{d}\Bigr)\cdot |\ln \e_d^2|}} =
    \dfrac{\sqrt{|\ln \e_d^2|}}{\sqrt{d}\cdot\sqrt{|\ln \omega|\cdot(d + c_\omega|\ln \e_d^2|)}}\\
    &&{}\quad \leqslant \dfrac{\sqrt{|\ln \e_d^2|}}{\sqrt{d}\cdot\sqrt{|\ln \omega|\cdot c_\omega|\ln \e_d^2|}}=\frac{1}{\sqrt{d}} \to 0,\quad d \to \infty.
\end{eqnarray*}
We return to $H_d$:
\begin{eqnarray*}
	 H_d(t)=\Bigl(1-  \tfrac{c_\omega a_d}{d}\, \bigl(   \exp\bigl\{\tfrac{it}{b_d}|\ln \omega|\bigr\}-1\bigr)\Bigr)^{-d} \Bigl(\exp\bigl\{it\tfrac{a_d}{db_d}\bigr\}\Bigr)^{-d},\quad t\in\R,\quad d\in\N.
\end{eqnarray*}
We fix $t\in\R$. Due to comments above, the following expansions are valid:
\begin{eqnarray*}
	\exp\bigl\{\tfrac{ita_d}{db_d}\bigr\} = 1 + it\tfrac{a_d}{db_d} - \tfrac{t^2}{2}\cdot \tfrac{a_d^2}{d^2 b_d^2} +
	O\Bigl(\tfrac{a_d^3}{d^3b_d^3}\Bigr),\quad d \to \infty,
\end{eqnarray*}
and
\begin{eqnarray*}
	\exp\bigl\{\tfrac{it}{b_d}|\ln \omega|\bigr\}-1= \tfrac{it}{b_d}|\ln \omega| -
	\tfrac{t^2}{2b_d^2}|\ln \omega|^2 + O\Bigl(\tfrac{1}{b_d^3}\Bigr), \quad d \to \infty.
\end{eqnarray*}
From the latter we get
\begin{eqnarray*}
	1-  \tfrac{c_\omega a_d}{d}\, \bigl(   \exp\bigl\{\tfrac{it}{b_d}|\ln \omega|\bigr\}-1\bigr)&=&	1-  \tfrac{c_\omega a_d}{d}\, \Bigl( \tfrac{it}{b_d}|\ln \omega| -
	\tfrac{t^2}{2b_d^2}|\ln \omega|^2 + O\Bigl(\tfrac{1}{b_d^3}\Bigr) \Bigr)\\
	&=& 1- it \tfrac{a_d}{d b_d}+ \tfrac{t^2}{2}\cdot\tfrac{a_d}{d b_d^2}\,|\ln \omega|+  O\Bigl(\tfrac{a_d}{d b_d^3}\Bigr),\quad d \to \infty.
\end{eqnarray*}
Thus we obtain
\begin{eqnarray*}
	H_d(t)=\Biggr[\biggl(  1- it \tfrac{a_d}{d b_d}+ \tfrac{t^2}{2}\cdot\tfrac{a_d}{d b_d^2}\,|\ln \omega|+  O\Bigl(\tfrac{a_d}{d b_d^3}\Bigr)\biggr)\cdot \biggl(1 + it\tfrac{a_d}{db_d} - \tfrac{t^2}{2}\cdot \tfrac{a_d^2}{d^2 b_d^2} +
	O\Bigl(\tfrac{a_d^3}{d^3b_d^3}\Bigr)\biggr)\Biggr]^{-d},\quad d\to\infty.
\end{eqnarray*}
We next write
\begin{eqnarray*}
	H_d(t)&=&\biggl[1- it \tfrac{a_d}{d b_d}+ it \tfrac{a_d}{d b_d} + \tfrac{t^2}{2}\cdot\tfrac{a_d}{d b_d^2}\,|\ln \omega|- \tfrac{t^2}{2}\cdot \tfrac{a_d^2}{d^2 b_d^2}-\bigl( it \tfrac{a_d}{d b_d} \bigr)^2
	+\tfrac{t^2}{2}\cdot\tfrac{a_d}{d b_d^2}\,|\ln \omega|\cdot  it\tfrac{a_d}{db_d} + it \tfrac{a_d}{d b_d}\cdot \tfrac{t^2}{2}\cdot \tfrac{a_d^2}{d^2 b_d^2}\\&&{}\quad - \tfrac{t^2}{2}\cdot\tfrac{a_d}{d b_d^2}\,|\ln \omega|\cdot \tfrac{t^2}{2}\cdot \tfrac{a_d^2}{d^2 b_d^2}+O\Bigl(\tfrac{a_d^3}{d^3b_d^3}\Bigr) +O\Bigl(\tfrac{a_d}{d b_d^3}\Bigr)\biggr]^{-d}\\
	&=& \biggl[1 + \tfrac{t^2}{2}\cdot\Bigl(\tfrac{a_d}{d b_d^2}\,|\ln \omega|+  \tfrac{a_d^2}{d^2 b_d^2}\Bigr) +  \tfrac{it^3}{2}\cdot\Bigl(\tfrac{a_d^2}{d^2 b_d^3}\,|\ln \omega|+  \tfrac{a_d^3}{d^3 b_d^3}\Bigr)\\
	&&{}\quad - \tfrac{t^4}{4}\cdot\tfrac{a_d^3}{d^3 b_d^4}\,|\ln \omega|+O\Bigl(\tfrac{a_d^3}{d^3b_d^3}\Bigr) +O\Bigl(\tfrac{a_d}{d b_d^3}\Bigr)\biggr]^{-d},\quad d\to\infty.
\end{eqnarray*}
Here
\begin{eqnarray*}
	\dfrac{a_d}{d b_d^2}\,|\ln \omega|+  \dfrac{a_d^2}{d^2 b_d^2}&=&|\ln \omega|\cdot\dfrac{a_d^2}{d^2 b_d^2}\cdot \dfrac{d+c_\omega a_d}{a_d}\\
	&=&|\ln \omega|\cdot\dfrac{|\ln \e_d^2|}{d\cdot|\ln \omega|\cdot(d + c_\omega|\ln \e_d^2|)}\cdot \dfrac{d+c_\omega|\ln \e_d^2|}{|\ln \e_d^2|}=\dfrac{1}{d},\quad d\in\N.
\end{eqnarray*}
Since 
\begin{eqnarray*}
	\dfrac{a_d}{db_d}=O(d^{-1/2}),\quad d\to\infty,
\end{eqnarray*}
we have
\begin{eqnarray*}
	\tfrac{it^3}{2}\cdot\Bigl(\tfrac{a_d^2}{d^2 b_d^3}\,|\ln \omega|+  \tfrac{a_d^3}{d^3 b_d^3}\Bigr)=\tfrac{it^3}{2}\cdot\Bigl(O\Bigl( \tfrac{1}{d b_d}\Bigr)+O\Bigl( \tfrac{1}{d\sqrt{d}}\Bigr)\Bigr)=o\bigl(\tfrac{1}{d}\bigr),\quad d\to\infty,
\end{eqnarray*}
and
\begin{eqnarray*}
	-\tfrac{t^4}{4}\cdot\dfrac{a_d^3}{d^3 b_d^4}\,|\ln \omega|+O\Bigl(\tfrac{a_d^3}{d^3b_d^3}\Bigr)= O\Bigl( \tfrac{1}{d\sqrt{d} b_d}\Bigr)+O\Bigl( \tfrac{1}{d\sqrt{d}}\Bigr)=o\bigl(\tfrac{1}{d}\bigr),\quad d\to\infty.
\end{eqnarray*}
Observe that
\begin{eqnarray*}
    \dfrac{a_d}{db_d^3} = \dfrac{|\ln \e_d^2|}{d|\ln \omega|^{3/2}\cdot\Bigl(1+ \tfrac{c_\omega|\ln \e_d^2|}{d}\Bigr)^{3/2}\cdot |\ln \e_d^2|^{3/2}}
    \leqslant \dfrac{1}{d\cdot|\ln \omega|^{3/2}\cdot  |\ln \e_d^2|^{1/2}}=o\bigl(\tfrac{1}{d}\bigr),\quad d\to\infty.
\end{eqnarray*}
Thus
\begin{eqnarray*}
    H_d(t) = \left(1 + \tfrac{t^2}{2d} + o\bigl(\tfrac{1}{d}\bigr)\right)^{-d} \to e^{-t^2 / 2},\quad d \to \infty.
\end{eqnarray*}
Recall that $\exp\{-t^2/2\}$ is the characteristic function corresponding to the standard normal law. By the well-known L\'evy's continuity theorem (see \cite{Lukacs} p. 48--49), we have the pointwise convergence
\begin{eqnarray}\label{conv_PPhi}
	\Probab\Biggl(\dfrac{S_{d, \e_d} - \Expec S_{d, \e_d}}{\sqrt{\Var S_{d, \e_d}}}\leqslant x\Biggr) \to\Phi(x), \quad d \to \infty, 
\end{eqnarray}
for any $x$ from the set $\mathcal{C}(\Phi)$ of all continuity points of $\Phi$. Since $\Phi$ is a  continuous function on $\R$, the set $\mathcal{C}(\Phi)$ is exactly $\R$ and, moreover, the convergence \eqref{conv_PPhi} is uniform on $\R$ (see \cite{Petrov}, p. 11).  Thus we come to \eqref{normal_proposition_conv}.\quad $\Box$\\

\section{Proofs of the main results}

\textbf{Proof of Theorem \ref{fixed_error_theorem}}.\quad We fix $\e\in(0,1)$. Let us use formula \eqref{eq_fourth_step} and also \eqref{const-omega-eqnarray}:
\begin{eqnarray*}
	n(d, \e) &=& \Bigl(1+\tfrac{c_\omega|\ln \e^2|}{d} \Bigr)^d\,
	\Expec \biggl[\exp\Big\{c_\omega\ln\Bigl(1+\tfrac{d}{c_\omega|\ln \e^2|}\Bigr)S_{d, \e}\Big\}
	\id\Big(S_{d, \e} < |\ln \e^2|\Big) \biggr]\\
	 &=&\Bigl(1+\tfrac{c_\omega|\ln \e^2|}{d} \Bigr)^d\sum_{\substack{m \in \N_0:\\ m|\ln\omega| <|\ln \e^2|}}
	\exp\Big\{c_\omega\ln\Bigl(1+\tfrac{d}{c_\omega|\ln \e^2|}\Bigr)m|\ln \omega|\Big\}\cdot\Probab\bigl(S_{d, \e} = m|\ln \omega|\bigr)\\
	&=&\Bigl(1+\tfrac{c_\omega|\ln \e^2|}{d} \Bigr)^d\sum_{\substack{m \in \N_0:\\ m < c_\omega|\ln \e^2|}}
	\exp\Big\{\ln\Bigl(1+\tfrac{d}{c_\omega|\ln \e^2|}\Bigr)m\Big\}\cdot\Probab\bigl(S_{d, \e} = m|\ln \omega|\bigr)\\
	&=&\Bigl(1+\tfrac{c_\omega|\ln \e^2|}{d} \Bigr)^d\sum_{m=0}^{N_\omega(\e)}
	\Bigl(1+\tfrac{d}{c_\omega|\ln \e^2|}\Bigr)^m\Probab\bigl(S_{d, \e} = m|\ln \omega|\bigr),
\end{eqnarray*}
where we set $N_{\omega}(\e)\defeq\max\bigl\{m\in\N_0: m< c_\omega |\ln \e^2|\bigr\}$ as in the statement of the theorem. Due to \eqref{conc_1plusedd}, we have
\begin{eqnarray*}
	\Bigl(1+\tfrac{c_\omega|\ln \e^2|}{d} \Bigr)^d=e^{c_\omega |\ln\e^2|}\cdot \Bigl(1+O\bigl(\tfrac{1}{d}\bigr)\Bigr),\quad d\to\infty.
\end{eqnarray*}
Applying Proposition \ref{prop-poisson} to the sum above, we obtain
\begin{eqnarray*}
	\sum_{m=0}^{N_\omega(\e)}
	\Bigl(1+\tfrac{d}{c_\omega|\ln \e^2|}\Bigr)^m\Probab\bigl(S_{d, \e} = m|\ln \omega|\bigr) &=& \sum_{m=0}^{N_\omega(\e)}\biggl[
	\Bigl(1+\tfrac{d}{c_\omega|\ln \e^2|}\Bigr)^m\Bigl(e^{-c_\omega|\ln \e^2|}\cdot\tfrac{(c_\omega|\ln \e^2|)^m}{m!}+O\bigl(\tfrac{1}{d}\bigr)\Bigr)\biggr]\\
	&=& \Bigl(1+\tfrac{d}{c_\omega|\ln \e^2|}\Bigr)^{N_{\omega}(\e)} e^{-c_\omega|\ln \e^2|}\cdot\dfrac{(c_\omega|\ln \e^2|)^{N_{\omega}(\e)}}{N_{\omega}(\e)!}+ O\bigl(d^{N_\omega(\e)-1}\bigr)\\
	&=& \Bigl(\tfrac{d}{c_\omega|\ln \e^2|}\Bigr)^{N_{\omega}(\e)} e^{-c_\omega|\ln \e^2|}\cdot\dfrac{(c_\omega|\ln \e^2|)^{N_{\omega}(\e)}}{N_{\omega}(\e)!}+ O\bigl(d^{N_\omega(\e)-1}\bigr)\\
	&=&  e^{-c_\omega|\ln \e^2|}\cdot\dfrac{d^{N_{\omega}(\e)}}{N_{\omega}(\e)!}+ O\bigl(d^{N_\omega(\e)-1}\bigr),\quad d\to\infty.
\end{eqnarray*}
Thus
\begin{eqnarray*}
	n(d,\e)&=& e^{c_\omega |\ln\e^2|}\cdot \Bigl(1+O\bigl(\tfrac{1}{d}\bigr)\Bigr)\cdot\biggl( e^{-c_\omega|\ln \e^2|}\cdot\dfrac{d^{N_{\omega}(\e)}}{N_{\omega}(\e)!}+ O\bigl(d^{N_\omega(\e)-1}\bigr)\biggr)\\
	&=& \dfrac{d^{N_{\omega}(\e)}}{N_{\omega}(\e)!}+ O\bigl(d^{N_\omega(\e)-1}\bigr),\quad d\to\infty,
\end{eqnarray*}
as required.\quad $\Box$\\

\textbf{Proof of Theorem \ref{logarithmic_nonfixed_error_theorem}.}\quad  Suppose that the sequence  $(\e_d)_{d\in\N}$ is given.  Due to \eqref{eq_fourth_step} we have
\begin{eqnarray*}
    n(d, \e_d) = \Bigl(1+\tfrac{c_\omega|\ln \e_d^2|}{d} \Bigr)^d\,
	\Expec \biggl[\exp\Big\{c_\omega\ln\Bigl(1+\tfrac{d}{c_\omega|\ln \e_d^2|}\Bigr)S_{d, \e_d}\Big\}
	\id\Big(S_{d, \e_d} < |\ln \e^2_d|\Big)\biggr],\quad d\in\N.
\end{eqnarray*}
From this we easily obtain the upper bound for $n(d, \e_d)$:
\begin{eqnarray*}
	n(d, \e_d) &\leqslant& \Bigl(1+\tfrac{c_\omega|\ln \e_d^2|}{d} \Bigr)^d\,\exp\Big\{c_\omega\ln\Bigl(1+\tfrac{d}{c_\omega|\ln \e_d^2|}\Bigr)|\ln \e^2_d|\Bigl\}
	\Expec 
	\id\bigl(S_{d, \e_d} < |\ln \e^2_d|\bigr)\\
	&=& \exp\biggl\{d\cdot\ln\Bigl(1+\tfrac{c_\omega|\ln \e_d^2|}{d}\Bigr) +
	c_\omega|\ln \e_d^2|\cdot\ln\Bigl(1+\tfrac{d}{c_\omega|\ln \e_d^2|}\Bigr)\biggr\} \Probab\bigl(S_{d, \e_d} < |\ln \e^2_d|\bigr)\\
	&\leqslant&\exp\biggl\{d\cdot\ln\Bigl(1+\tfrac{c_\omega|\ln \e_d^2|}{d}\Bigr) +
	c_\omega|\ln \e_d^2|\cdot\ln\Bigl(1+\tfrac{d}{c_\omega|\ln \e_d^2|}\Bigr)\biggr\}, \quad d\in\N.
\end{eqnarray*}
For the lower bound we introduce $b_d\defeq \sqrt{\Var S_{d,\e_d}}$, $d\in\N$, as in \eqref{def_adbd}. So we have
\begin{eqnarray*}
	n(d, \e_d) &\geqslant& \Bigl(1+\tfrac{c_\omega|\ln \e_d^2|}{d} \Bigr)^d\,
	\Expec \biggl[\exp\Big\{c_\omega\ln\Bigl(1+\tfrac{d}{c_\omega|\ln \e_d^2|}\Bigr)S_{d, \e_d}\Big\}
	\id\Big(|\ln \e^2_d|-b_d<S_{d, \e_d} < |\ln \e^2_d|\Big)\biggr]\\
	&\geqslant& \Bigl(1+\tfrac{c_\omega|\ln \e_d^2|}{d} \Bigr)^d\,\exp\Big\{c_\omega\ln\Bigl(1+\tfrac{d}{c_\omega|\ln \e_d^2|}\Bigr)\bigl(|\ln \e^2_d|-b_d\bigr)\Big\}
	\Expec \id\Big(|\ln \e^2_d|-b_d<S_{d, \e_d} < |\ln \e^2_d|\Big)\\
	&=& \exp\biggl\{d\cdot\ln\Bigl(1+\tfrac{c_\omega|\ln \e_d^2|}{d}\Bigr) +
	c_\omega|\ln \e_d^2|\cdot\ln\Bigl(1+\tfrac{d}{c_\omega|\ln \e_d^2|}\Bigr)\biggr\}\\
	&&{}\cdot \exp\Big\{-c_\omega b_d\ln\Bigl(1+\tfrac{d}{c_\omega|\ln \e_d^2|}\Bigr)\Big\}\Probab\Big(|\ln \e^2_d|-b_d<S_{d, \e_d} < |\ln \e^2_d|\Big),\quad d\in\N.
\end{eqnarray*}
According to \eqref{conc_adbd}, observe that
\begin{eqnarray*}
	\dfrac{b_{d}}{|\ln \e_d^2|} = \dfrac{1}{|\ln \e_d^2|}\sqrt{|\ln \omega|\cdot\Bigl(1+ \tfrac{c_\omega|\ln \e_d^2|}{d}\Bigr)\cdot |\ln \e_d^2|}=\sqrt{|\ln \omega|\cdot\Bigl(\tfrac{1}{|\ln \e_d^2|}+ \tfrac{c_\omega}{d}\Bigr)}
	\to 0, \quad d \to \infty.
\end{eqnarray*}
Therefore
\begin{eqnarray*}
	-c_\omega b_d\ln\Bigl(1+\tfrac{d}{c_\omega|\ln \e_d^2|}\Bigr)=-\dfrac{b_{d}}{|\ln \e_d^2|}\cdot c_\omega |\ln \e_d^2|\cdot \ln\Bigl(1+\tfrac{d}{c_\omega|\ln \e_d^2|}\Bigr)=o\Bigl(|\ln \e_d^2|\cdot \ln\Bigl(1+\tfrac{d}{c_\omega|\ln \e_d^2|}\Bigr)\Bigr), \quad d \to \infty.
\end{eqnarray*}
Next, since $|\ln \e_d^2|=\Expec S_{d, \e_d}$ and $b_d=\sqrt{\Var S_{d,\e_d}}$,  by Proposition \ref{normal_proposition}, we have
\begin{eqnarray*}
    \Probab\Big(|\ln \e^2_d|-b_d<S_{d, \e_d} < |\ln \e^2_d|\Big)=
    \Probab\Biggl(-1<\dfrac{S_{d, \e_d} - \Expec S_{d, \e_d}}{\sqrt{\Var S_{d, \e_d}}}<0\Biggr)
    =\Phi(0) - \Phi(-1)+o(1),\quad d\to\infty.
\end{eqnarray*}
Here $C\defeq \Phi(0)-\Phi(-1)$ is a positive constant:
\begin{eqnarray*}
	C=\dfrac{1}{\sqrt{2\pi }}\int_{-1}^{0} e^{-x^2/2}\, \dd x=\dfrac{1}{\sqrt{2\pi }}\int_{0}^{1} e^{-x^2/2}\, \dd x>0.
\end{eqnarray*}
In particular, this means that
\begin{eqnarray*}
	\ln\Probab\Big(|\ln \e^2_d|-b_d<S_{d, \e_d} < |\ln \e^2_d|\Big)=\ln\bigl(C+o(1)\bigr)=o\Bigl(|\ln \e_d^2|\cdot \ln\Bigl(1+\tfrac{d}{c_\omega|\ln \e_d^2|}\Bigr)\Bigr),\quad d\to\infty.
\end{eqnarray*}
Thus the lower bound admits the form
\begin{eqnarray*}
	n(d,\e_d)\geqslant\exp\biggl\{d\cdot\ln\Bigl(1+\tfrac{c_\omega|\ln \e_d^2|}{d}\Bigr) +
	c_\omega|\ln \e_d^2|\cdot\ln\Bigl(1+\tfrac{d}{c_\omega|\ln \e_d^2|}\Bigr)+o\Bigl(|\ln \e_d^2|\cdot \ln\Bigl(1+\tfrac{d}{c_\omega|\ln \e_d^2|}\Bigr)\Bigr)\biggr\}.
\end{eqnarray*}

The obtained bounds imply the required result. \quad $\Box$

\end{document}